\documentclass{amsart}

\usepackage{amssymb}
\usepackage{amscd}

\def\frak{\mathfrak}
\def\Bbb{\mathbb}

\newtheorem{thm}{Theorem}
\newtheorem*{thm*}{Theorem}
\newtheorem{prop}[thm]{Proposition}
\newtheorem*{prop*}{Proposition}
\newtheorem{lem}[thm]{Lemma}
\newtheorem*{lem*}{Lemma}

\newtheorem*{kor*}{Corollary}

\renewcommand{\exp}{\operatorname{exp}}
\newcommand{\id}{\operatorname{id}}

\newcommand{\End}{\operatorname{End}}

\newcommand{\Ups}{\Upsilon}
\newcommand{\Rho}{\mathsf{P}}

\newcommand{\fg}{{\frak g}}
\newcommand{\x}{\times}

\let\ccdot\cdot
\def\cdot{\hbox to 2.5pt{\hss$\ccdot$\hss}}

\newcommand{\al}{\alpha}
\newcommand{\be}{\beta}
\newcommand{\ga}{\gamma}

\newcommand{\la}{\lambda}

\renewcommand{\phi}{\varphi}

\newcommand{\Ga}{\Gamma}
\newcommand{\La}{\Lambda}
\newcommand{\Om}{\Omega}

\begin{document}
\title{AHS--structures and affine holonomies}

\author{Andreas \v Cap}

\address{Institut f\"ur Mathematik, Universit\"at Wien, Nordbergstra\ss
  e 15, A--1090 Wien, Austria and International Erwin Schr\"odinger
  Institute for Mathematical Physics, Boltzmanngasse 9, A--1090 Wien,
  Austria} 

\email{Andreas.Cap@esi.ac.at} 

\thanks{Supported by
    project P 19500--N13 of the ``Fonds zur F\"orderung der
    wissenschaftlichen Forschung'' (FWF)}



\commby{J. Wolfson}


\begin{abstract} 
  We show that a large class of non--metric, non--symplectic affine
  holonomies can be realized, uniformly and without case by case
  considerations, by Weyl connections associated to the natural
  AHS--structures on certain generalized flag manifolds.
\end{abstract}

\maketitle

\section{Introduction}\label{1}

The classification of the possible irreducible holonomies of
non--locally symmetric torsion free affine connections is a
cornerstone in differential geometry. A list of possible holonomy Lie
algebras was compiled by M.~Berger, see \cite{Berger}, later a few
small corrections and several extensions to this list were found. It
took quite a long time until existence of all these holonomies was
proved, often case by case and uncovering interesting relations to
several areas of differential geometry and related fields. The program
was finally completed by S.~Merkulov and L.~Schwachh\"ofer in
\cite{Merkulov-Schwachhoefer}.

Apart from the case of metric holonomies (i.e.~holonomies of
connections admitting a parallel pseudo--Riemannian metric), the
classification of possible affine holonomies has surprisingly close
relations to the classification of certain types of symmetric spaces
or, equivalently, of certain types of parabolic subalgebras in simple
Lie algebras respectively certain generalized flag manifolds. In the
original existence proofs for holonomies, these relations were not
exploited systematically. For the case of symplectic holonomies
(i.e.~holonomies of connections admitting a parallel symplectic form)
this was done later by M.~Cahen and L.~Schwachh\"ofer in
\cite{Cahen-Schwachhoefer}. There the authors construct special
symplectic connections starting from certain generalized flag
manifolds. This not only provides examples of all symplectic
holonomies, but locally produces all connections with such holonomies.

The aim of this article is to give a conceptual proof of existence of
affine ho\-lo\-no\-mies, which exploits the relation to parabolic
subalgebras, in an easier case, namely for the holonomies related to
the classification of Hermitian symmetric spaces. Here the simplest
description is in terms of so--called $|1|$--graded simple Lie
algebras. For a $|1|$--grading $\fg=\fg_{-1}\oplus\fg_0\oplus\fg_1$ of
a simple Lie algebra $\fg$, the subspace $\fg_0$ turns out to be
reductive with one--dimensional center. The adjoint action defines a
representation of $\fg_0$ on the vector space $\fg_{-1}$. Our main
result is:
\begin{thm}\label{main}
  For any $|1|$--graded simple Lie algebra
  $\fg=\fg_{-1}\oplus\fg_0\oplus\fg_1$, the adjoint representations of
  $\fg_0$ and of its semisimple part on $\fg_{-1}$ can be realized as
  holonomy Lie algebras of torsion free, non locally symmetric affine
  connections on a compact manifold.
\end{thm}
The classification of $|1|$--gradings on simple Lie algebras and hence
of the ho\-lo\-no\-mies covered by this construction is well known and can
be found in table 1 below. Note that in contrast to the results of
\cite{Cahen-Schwachhoefer} for symplectic holonomies, which are local
in nature, we obtain global connections on compact manifolds, indeed
on generalized flag manifolds.

A Lie group with $|1|$--graded simple Lie algebra determines a
geometric structure, namely a first order structure with a certain
structure group $G_0$, which can be canonically prolonged to a normal
Cartan geometry. These are the so--called AHS--structures, with
conformal structures providing the prototypical example. Any structure
of this type (as well as the more general parabolic geometries) comes
with a class of distinguished affine connections, called Weyl
connections. The theory of Weyl structures (which are equivalent
descriptions of the Weyl connections) as developed in \cite{Weyl}
provides the main technical input for our results.

\subsection*{Acknowledgements} The idea to use Weyl connections for
AHS structures to realize special holonomies was originally brought up
several years ago by J.~Slov\'ak and D.~Alekseevsky. Conversations
with S.~Armstrong and J.~Slov\'ak have been very helpful.

\section{AHS--structures and Weyl connections}\label{2}

A $|1|$--grading on a simple Lie algebra $\fg$ is a decomposition
$\fg=\fg_{-1}\oplus\fg_0\oplus\fg_1$ which defines a grading of $\fg$,
i.e.~is such that $[\fg_i,\fg_j]\subset\fg_{i+j}$ where $\fg_k=0$ for
$k\notin\{-1,0,1\}$. We will always assume that the rank of $\fg$ is
bigger than one, so $\dim(\fg_{-1})>1$. It turns out that the Lie
subalgebra $\fg_0\subset\fg$ is always reductive with one--dimensional
center. The semisimple part of $\fg_0$ will be denoted by
$\fg_0^{ss}$. By the grading property, the adjoint representation of
$\fg$ restricts to representations of $\fg_0$ on $\fg_{-1}$ and
$\fg_1$. These representations are always irreducible and dual to each
other via the Killing form of $\fg$. Finally, $\frak
p:=\fg_0\oplus\fg_1$ is a maximal parabolic subalgebra of $\fg$ with
nilradical $\fg_1$. 

$|1|$--gradings on simple Lie algebras are closely related to
Hermitian and para--Hermitian symmetric spaces, as well as to
questions of Lie algebras with non--trivial prolongations. The full
classification of is well known, see e.g.~\cite{Kobayashi-Nagano}.
\begin{prop}\label{prop1}
  Table 1 lists all real and complex simple $|1|$--graded Lie algebras
  $\fg$ of rank bigger than one, together with the subalgebras
  $\fg_0$ and the representations on $\fg_{-1}$, with $\Bbb K$
  denoting $\Bbb R$ or $\Bbb C$.
\end{prop}

\begin{table}[t]
  \centering
  \begin{tabular}{|c|c|c|}
\hline
$\fg$ & $\fg_0$ & $\fg_{-1}$\\
\hline
\hline
$\frak{sl}(n+1,\Bbb K)$, $n\geq 2$ & $\frak{gl}(n,\Bbb K)$ & $\Bbb K^n$\\
\hline
$\frak{sl}(p+q,\Bbb K)$, $p,q\geq 2$ & $\frak s(\frak{gl}(p,\Bbb
K)\x\frak{gl}(q,\Bbb K))$ & $\Bbb K^{p*}\otimes\Bbb K^q$ \\ 
\hline 
$\frak{sl}(p+q,\Bbb H)$, $p,q\geq 1$ & $\frak s_{\Bbb
  R}(\frak{gl}(p,\Bbb H)\x\frak{gl}(q,\Bbb H))$ & 
$L_{\Bbb H}(\Bbb H^p,\Bbb H^q)$\\
\hline
$\frak{su}(p,p)$, $p\geq 2$ & $\frak s_{\Bbb R}(\frak{gl}(p,\Bbb C))$ & $\frak
u(p)$\\ 
\hline
$\frak{sp}(2n,\Bbb K)$, $n\geq 3$ & $\frak{gl}(n,\Bbb K)$ & $S^2\Bbb K^n$ \\ 
\hline
$\frak{so}(p+1,q+1)$, $p+q\geq 3$ & $\frak{cso}(p,q)$ & $\Bbb R^{p+q}$ \\
\hline
$\frak{so}(n+2,\Bbb C)$, $n\geq 3$ & $\frak{cso}(n,\Bbb C)$ & $\Bbb C^n$ \\ 
\hline
$\frak{so}(n,n)$, $n\geq 4$ & $\frak{gl}(n,\Bbb R)$ & $\La^2\Bbb
R^n$\\ 
\hline
$\frak{so}(2n,\Bbb C)$, $n\geq 4$ & $\frak {gl}(n,\Bbb C)$ & $\La^2\Bbb
C^n$\\ 
\hline
$\frak{so^*}(4n)$, $n\geq 2$ & $\frak{gl}(n,\Bbb H)$ & $\La^2\Bbb H^n$
\\
\hline
$EI$ (split $E_6$) & $\frak{cspin}(5,5)$ & $\Bbb R^{16}$ \\
\hline
$EIV$ (non--split $E_6$) & $\frak{cspin}(9,1)$ &   $\Bbb R^{16}$ \\
\hline
$E_6$ (complex) & $\frak{cspin}(10,\Bbb C)$ & $\Bbb C^{16}$ \\ 
\hline
$EV$ (split $E_7$) & $EI\oplus\Bbb R$ & $\Bbb R^{27}$ \\
\hline
$EVII$ (non--split $E_7$) &  $EIV\oplus\Bbb R$ & $\Bbb R^{27}$ \\
\hline
$E_7$ (complex) &  $E_6\oplus\Bbb C$ & $\Bbb C^{27}$ \\
\hline
  \end{tabular}
  \caption{Real and complex $|1|$--graded simple Lie algebras}
  \label{tab:1}
\end{table}

The algebras $\fg_0$ in Table 1 and their semisimple parts exhaust all
of Berger's original list of non--metric holonomies (Table 2 of
\cite{Merkulov-Schwachhoefer}) except the full symplectic algebras, as
well as some exotic holonomies (Table 3 of
\cite{Merkulov-Schwachhoefer}).

\medskip

Suppose that we have given a $|1|$--graded simple Lie algebra $\fg$
and a Lie group $G$ with Lie algebra $\fg$. Then there are natural
subgroups $G_0\subset P\subset G$ with Lie algebras $\fg_0$ and $\frak
p$. Namely, we let $P$ consist of all elements whose adjoint actions
preserve $\frak p$ and $\fg_1$, while the adjoint actions of elements
of $G_0$ preserve any of the grading components $\fg_i$. It turns our
that the exponential mapping restricts to a diffeomorphism from
$\fg_1$ onto a closed normal subgroup $P_+\subset P$ and that $P$ is
the semidirect product of $G_0$ and $P_+$. 

The representation $G_0\to GL(\fg_{-1})$ obtained from the adjoint
representation is infinitesimally injective, so the notion of a first
order structure with structure group $G_0$ makes sense on manifolds of
dimension $\dim(\fg_{-1})$.  The structures obtained in that way are
called \textit{AHS--structures}, \textit{generalized conformal
  structures}, \textit{irreducible parabolic geometries}, or
\textit{abelian parabolic geometries} in the literature. It turns out
that they are equivalent to normal Cartan geometries of type $(G,P)$,
see e.g.~\cite{AHS2}.

The \textit{homogeneous model} of the AHS--structure of type $(G,P)$
is the homogeneous space $G/P$. Here the first order $G_0$--structure
is given by the canonical projection $G/P_+\to G/P$, which is a
principal bundle with structure group $P/P_+=G_0$, and the soldering
form induced by the Maurer--Cartan form on $G$. The corresponding
Cartan geometry is the natural principal $P$--bundle $G\to G/P$ with
the Maurer--Cartan form as the Cartan connection. This Cartan
connection is flat by the Maurer--Cartan equation.  

\medskip

On any manifold endowed with an AHS--structure, there is a family of
preferred principal connections on the principal bundle defining the
$G_0$--structure. The theory of these connections is developed in
\cite{Weyl}. Let us briefly describe the case of conformal structures,
which motivated the whole theory. 

Let $M$ be a smooth manifold of dimension $\geq 3$ endowed with a
conformal class $[\ga]$ of Riemannian metrics. A \textit{Weyl
  connection} on $TM$ is a torsion free linear connection $\nabla$
such that $\nabla_\xi\ga=f\ga$ for one (or equivalently any) metric
$\ga$ in the conformal class and any vector field $\xi\in\frak X(M)$.
Here $f$ is some smooth function (depending on $\xi$). The
Levi--Civita connection of any of the metrics in the conformal class
is an example of a Weyl connection, but these do not exhaust all Weyl
connections. The Weyl connections can equivalently be considered as
principal connections on the conformal frame bundle.

One can consider the conformal class of metrics as a ray subbundle in
$S^2T^*M$, which can be viewed as a principal bundle with structure
group $\Bbb R_+$. This is called the bundle of scales. Any Weyl
connection induces a principal connection on the bundle of scales, and
it turns out that this induces a bijection between the set of Weyl
connections and the set of all principal connections on the bundle of
scales. The Levi--Civita connections of metrics in the conformal class
exactly correspond to the flat connections induced by global sections
of the bundle of scales. The space of principal connections on the
bundle of scales is an affine space modelled on the vector space
$\Om^1(M)$ of one--forms on $M$, so one can carry over the the affine
structure to the space of all Weyl connections (obtaining a
non--trivial action of one--forms on Weyl connections). The subspace
of Levi--Civita connections thereby becomes an affine space modelled
on the space of exact one--forms. The totally tracefree part of the
curvature is the same for all Weyl connections. This is the Weyl
curvature of the conformal structure. The trace part of the curvature
of a Weyl connection is best described by the \textit{Rho--tensor} (a
trace adjusted version of the Ricci--curvature). 

Now all this generalizes to all AHS--structures (and further to
parabolic geometries). One always has a distinguished family of
principal connections on the bundle defining the $G_0$--structure,
which is in bijective correspondence with the space of all connections
on a principal $\Bbb R_+$--bundle, and hence forms an affine space
modelled on one--forms. These are called \textit{Weyl structures} or
\textit{Weyl connections}. The Weyl structures coming from global
sections of the $\Bbb R_+$--bundle are called \textit{exact}. It will
be important in the sequel that (as for conformal structures) exact
Weyl connections actually are induced from a further reduction of the
principal bundle defining the $G_0$--structure. The structure group of
this reduced bundle has Lie algebra $\fg^{ss}_0$, the semisimple part
of $\fg_0$.  Exact Weyl structure form an affine space modelled on the
space of exact one--forms.

In the case of AHS--structures, the set of all Weyl connections is easy
to describe. There is a basic invariant of such a geometry called the
harmonic torsion, and the Weyl connections are exactly those
connections on the principal bundle defining the $G_0$--structure
which have that torsion. (For some geometries, like conformal
structures, this torsion always has to vanish.) Finally, any Weyl
structure comes with a Rho--tensor, a $T^*M$--valued one--form which
describes a part of the curvature of the Weyl connection. For our
purposes, the following following special case of these facts will be
sufficient. 
\begin{prop}\label{prop2}
  For the homogeneous model $G/P$ of an AHS structure, the Weyl
  connections are exactly the torsion free linear connections on
  $T(G/P)$ which are induced from the principal $G_0$--bundle
  $G/P_+\to G/P$. Any Weyl connection has holonomy Lie algebra
  contained in $\fg_0$, and for exact Weyl connections the holonomy
  Lie algebra is even contained in $\fg_0^{ss}$. 
\end{prop}

An important feature of the theory of Weyl structures is that there is
an explicit description of the behavior of all relevant quantities
under a change of Weyl structure, which is valid for all the geometries
in question. This needs some tensorial maps coming from the AHS
structure.

For a $G_0$--structure $E\to M$ one by definition has $TM\cong
E\x_{G_0}\fg_{-1}$. We have already noted that $\fg_1$ is dual to
$\fg_{-1}$ as a $G_0$--representation, so $T^*M\cong E\x_{G_0}\fg_1$.
Finally, $\fg_0$ can be viewed as a Lie subalgebra of
$L(\fg_{-1},\fg_{-1})$, so $E\x_{G_0}\fg_0$ can be naturally viewed as
a subbundle $\End_0(TM)$ of $L(TM,TM)$. Now the Lie bracket on $\fg$
is a $G_0$--equivariant map, so passing the associated bundles the
components of this bracket induce tensorial maps $TM\x
T^*M\to\End_0(TM)$, $\End_0(TM)\x TM\to TM$ and $\End_0(TM)\x T^*M\to
T^*M$, which we all denote by $\{\ ,\ \}$.

In terms of these brackets, one can now easily write a formula for the
affine structure on Weyl connections which is uniform for all the
AHS--structures as well as a formula for the change of the
Rho--tensor. The change of Weyl connections is most conveniently
expressed in terms of linear connections on associated vector bundles.
Given a representation $V$ of $G_0$ and the corresponding vector
bundle $F:=E\x_{G_0} V\to M$, the infinitesimal representation of
$\fg_0$ on $V$ induces a bilinear bundle map $\bullet :\End_0(TM)\x
F\to F$. The following formulae are taken from \cite[3.6]{Weyl}.
\begin{prop}\label{prop3}
  Let $E\to M$ be a first order $G_0$--structure, fix a Weyl structure
  and denote the induced connections on all associated bundles by
  $\nabla$. Let $\Ups\in\Om^1(M)$ be a one--form and let us indicate
  by hats the quantities associated to the Weyl structure obtained by
  modifying the initial structure by $\Ups$. Then for any vector field
  $\xi\in\frak X(M)$ we have

\noindent
(1) The modified Weyl connection on an associated vector bundle
$F=E\x_{G_0}V$ is given by
$$
\hat\nabla_\xi s=\nabla_\xi s-\{\Ups,\xi\}\bullet s,
$$
for all $s\in\Ga(F)$.

\noindent
(2) The Rho--tensor of the modified Weyl connection is given by
$$
\hat\Rho(\xi)=\Rho(\xi)+\nabla_\xi\Ups+\tfrac{1}{2}\{\Ups,\{\Ups,\xi\}\}. 
$$
\end{prop}

\section{Realizing affine holonomies}\label{3}

As we have seen in Proposition \ref{prop2}, for any Weyl connection
(respectively exact Weyl connection) on the homogeneous model $G/P$ of
an AHS--structure, the holonomy Lie algebra is contained in $\fg_0$
(respectively $\fg_0^{ss}$). Our aim is to show that there are Weyl
connections for which the holonomy Lie algebras equal these two
subalgebras. Our strategy will be to first construct an exact Weyl
connection on $G/P$ which is flat on an open neighborhood of $o=eP\in
G/P$. Then we exploit the affine structure on the spaces of Weyl
connections respectively exact Weyl connections.  We construct a one
form respectively an exact one form, such that modifying by these
one--forms one obtains connections with full holonomy Lie algebra. In
view of the explicit formulae for the changes of the data associated
to a Weyl structure under a change of Weyl structures, this can be
deduced from rather simple algebraic facts.

\begin{prop}\label{prop4}
  Let $G$ be a Lie group with $|1|$--graded simple Lie algebra $\fg$
  and let $G_0\subset P\subset G$ be the subgroups determined by the
  grading. Then there is a globally defined Weyl connection $\nabla$
  for $G/P$ which is flat with identically vanishing Rho--tensor
  locally around $o=eP\in G/P$.
\end{prop}
\begin{proof}
  By Proposition 3.2 of \cite{Weyl} there always exist global exact
  Weyl connections. On the other hand, since $\fg_{-1}$ is transversal
  to $\frak p$ in $\fg$, there is an open neighborhood $U$ of $0$ in
  $\fg_{-1}$ such that $X\mapsto \exp(X)g$ defines diffeomorphisms
  from $U\x P$ to an open neighborhood of $P$ in $G$ as well as from
  $U\x G_0$ to an open neighborhood of $G_0=P/P_+$ in $G/P_+$. Via
  these diffeomorphisms, the inclusion $G_0\hookrightarrow P$ defines
  a $G_0$--equivariant smooth section of $G\to G/P_+$, and hence a
  local Weyl structure as defined in \cite{Weyl}, over
  $V:=\exp(U)P\subset G/P$. (In fact, one may take $U=\fg_{-1}$, but
  this is not necessary for our purposes). The local trivialization of
  $G/P_+$ over $V$ gives rise to a local trivialization of any bundle
  of scales, and by construction the corresponding section of that
  bundle is parallel for the induced Weyl connection, so we have
  constructed an \textit{exact} local Weyl structure.
  
  Applying Proposition 2.10 of \cite{Cap-Schichl} to the opposite
  grading, we see that $\exp(\fg_{-1})G_0$ is a closed subgroup of $G$
  (the parabolic subgroup opposite to $P$). Our section thus is simply
  the restriction of the inclusion of this subgroup to some open
  subset, so the pull back of the Maurer--Cartan form of $G$ along
  this inclusion is the Maurer--Cartan form of the subgroup.  Hence it
  has values in $\fg_{-1}\oplus\fg_0$, which implies that the
  Rho--tensor vanishes identically, and the Maurer--Cartan equation
  shows that the resulting connection is flat.
  
  In view of the affine structure on exact Weyl connections, there is
  an exact one--form $df$ on $V$, which describes the change from the
  restriction of the initial global Weyl connection to our local flat
  Weyl connection. Now multiply $f$ by a smooth function with support
  in $V$ which is identically one locally around $o$. Then the result
  can be extended by zero to a globally defined smooth function on
  $G/P$, and the corresponding exact one--form transforms the initial
  global Weyl connection to an exact Weyl connection with the required
  properties. 
\end{proof}

  The classical example for this construction is again provided by
  conformal structures. Here $G/P$ is the sphere $S^n$, and one may
  use the round metric as the initial exact Weyl structure. The local
  flat Weyl structure is obtained by pulling back the flat metric on
  $\Bbb R^n$ via stereographic projection, which is well known to be
  conformal. Then one glues the two metrics using a partition of
  unity.

\medskip

Now let $\nabla$ be a Weyl connection as described in Proposition
\ref{prop4}, i.e.~the curvature $R$ and the Rho--tensor $\Rho$ of
$\nabla$ both vanish identically on an open neighborhood $U$ of
$o=eP\in G/P$. Now let us change the Weyl connection by a one--form
$\Ups\in\Om^1(G/P)$ and let us indicate the resulting quantities by
hats. Using that the geometry on $G/P$ is flat, Proposition 4.3 of
\cite{Weyl} shows that on $U$ we have
$$
\hat R(\xi,\eta)=(\partial\hat\Rho)(\xi,\eta)=
\{\hat\Rho(\xi),\eta\}-\{\hat\Rho(\eta),\xi\}. 
$$
On the other hand, $\hat\Rho$ is determined part (2) of Proposition
\ref{prop3}. Now the bundle map
$$
\partial:T^*(G/P)\otimes T^*(G/P)\to\La^2T^*(G/P)\otimes\End_0(T(G/P))
$$
is induced by a $G_0$--equivariant map on the representation spaces
inducing the bundles (which we will denote by the same letter), so it
is parallel for any Weyl connection. This implies that on $U$ we get
$\hat\nabla^k\hat R=(\id\otimes\partial)(\hat\nabla^k\hat\Rho)$. Now
we are ready to state the main algebraic input.

\begin{lem}\label{lemma1}
  The map $(\Rho,X,Y)\mapsto \partial\Rho(X,Y)$ induces surjections
  $S^2\fg_1\otimes\fg_{-1}\otimes\fg_{-1}\to \fg_0^{ss}$ and
  $\fg_1\otimes\fg_1\otimes\fg_{-1}\otimes\fg_{-1}\to \fg_0$.
\end{lem}
\begin{proof}
  Both maps are evidently $\fg_0$--equivariant, so it suffices to show
  that their images meet each irreducible component of the target
  space. Since $\dim(\fg_{-1})>1$, we can choose linearly independent
  elements $X,Y\in\fg_{-1}$ and a linear map $\Rho:\fg_{-1}\to\fg_1$
  such that $\Rho(Y)=0$ and $B(\Rho(X),Y)\neq 0$, where $B$ denotes
  the Killing form on $\fg$. Now it is well known that the center of
  $\fg_0$ is generated by the grading element $E$, i.e.~the element
  whose adjoint action on each $\fg_j$ is multiplication by $j$, and
  that the decomposition $\fg_0=\fg_0^{ss}\oplus\Bbb K E$ is orthogonal
  for $B$. Now by construction, $\partial\Rho(X,Y)=-[Y,\Rho(X)]$ and
  hence
$$
B(\partial\Rho(X,Y),E)=B(\Rho(X),[Y,E])=B(\Rho(X),Y)\neq 0.
$$
Hence the image of the second map is not contained in $\fg_0^{ss}$,
so it suffices to prove surjectivity of the first map.

Complexifying if necessary, we may assume that $\fg$ is a complex
$|1|$--graded simple Lie algebra, and then we can use the root
decomposition. (See \cite{Yamaguchi} for the description of
$|1|$--gradings in terms of roots.) There is a unique simple root
$\al$ such that the root space $\fg_\al$ is contained in $\fg_1$. Now
the Dynkin diagram of $\fg_0^{ss}$ is obtained by removing the node
corresponding to $\al$ and all edges connecting to this node in the
Dynkin diagram of $\fg$. Hence any simple ideal of $\fg_0^{ss}$
contains the root space $\frak g_{\be}$ for some simple root $\be$,
such that the nodes corresponding to $\al$ and $\be$ in the Dynkin
diagram of $\fg$ are connected. Equivalently, this means that $\be$ is
not orthogonal to $\al$ and hence $\al+\be$ is a root. By
construction, the root space $\fg_{\al+\be}$ is contained in $\fg_1$.
Now choose a basis $\{X_\ga\}$ of $\fg_{-1}$ such that each $X_\ga$
lies in the root space $\fg_{-\ga}$ and let $\{Z_\ga\}$ be the dual
basis of $\fg_1$.  Then for $\Rho= Z_\al^2+Z_{\al+\be}^2\in S^2\fg_1$ we get
$$
\partial\Rho(X_{\al},X_{\al+\be})=[X_\al,Z_{\al+\be}]-[X_{\al+\be},Z_\al].
$$
This is the sum of a nonzero element of $\fg_\be$ and a nonzero
element of $\fg_{-\be}$, so the image of our map meets the simple
ideal containing $\fg_{\be}$.
\end{proof}

\begin{lem}\label{lemma2}
  There is an element of $S^6\fg_{-1}^*$ which, interpreted as a
  linear map $S^4\fg_{-1}\to S^2\fg_{-1}^*$ is surjective. Likewise,
  there is an element in $S^5\fg_{-1}^*\otimes\fg_{-1}^*$ which is
  surjective when viewed as a map
  $S^4\fg_{-1}\to\fg_{-1}^*\otimes\fg_{-1}^*$. 
\end{lem}
\begin{proof}
  Let $\{e_i\}$ be a basis of $\fg_{-1}$ with dual basis $\{\la_i\}$,
  and consider $\sum_{i<j}\la_i^4\la_j^2\in S^6\fg_{-1}^*$. Then
  for $i<n=\dim(\fg_{-1})$, we obtain $\la_i^2$ as the image of
  $e_i^2e_{i+1}^2$, while $\la_n^2$ is the image of $e_{n-1}^4$.
  For $i<j$ we get $\la_i\la_j$ as the image of $e_i^3e_j$.
  
  On the other hand, consider $\sum_{i,j}\la_i^3\la_j^2\otimes
  \la_i\in S^5\fg_{-1}^*\otimes\fg_{-1}^*$. The corresponding map
  sends $e_i^4$ to $\la_i\otimes \la_i$ and $e_i^3e_j$ to
  $\la_j\otimes \la_i$.
\end{proof}

\begin{thm}\label{thm1}
  The space $G/P$ admits global Weyl connections with holonomy Lie
  algebra $\fg_0$ as well as global exact Weyl connections with
  holonomy Lie algebra $\fg_0^{ss}$. 
\end{thm}
\begin{proof}
  Let $\nabla$ be an exact Weyl connection on $G/P$ as in Proposition
  \ref{prop4}, and change it by $\Ups\in\Om^1(M)$. Then we have noted
  already that $\hat\nabla^k\hat
  R(o)=(\id\otimes\partial)(\hat\nabla^k\hat\Rho(o))$. Now suppose in
  addition that $\Ups$ has vanishing $k$--jet in $o$. Then part (2) of
  Proposition \ref{prop3} shows that $j^{k-1}_o\hat\Rho=0$ and that
  $\hat\nabla^k\hat\Rho(o)=\hat\nabla^k\nabla\Ups(o)$. Part (1) of
  Proposition \ref{prop3} inductively shows that
  $\hat\nabla^k\nabla\Ups(o)=\nabla^{k+1}\Ups (o)$. Of course we can
  find a global smooth one form $\Ups$ with $j^k_o\Ups=0$ such that
  $\nabla^{k+1}\Ups(o)$ (which is automatically symmetric in the first
  $k+1$ entries by flatness of $\nabla$) is any prescribed element of
  $S^{k+1}T_o(G/P)^*\otimes T_o(G/P)^*$. Likewise, this can be done
  for an exact one--form, with any prescribed element of
  $S^{k+2}T_o(G/P)^*$.
  
  Doing this for $k=4$ with elements which (via the isomorphism
  $T_o(G/P)\cong\fg_{-1}$) correspond to the ones described in Lemma
  \ref{lemma2}, Lemma \ref{lemma1} shows that we get a Weyl
  connection for which the values of $\hat\nabla^4\hat R$ fill all of
  $\fg_0$ and an exact Weyl connection for which these values fill all
  of $\fg_0^{ss}$. Since these values are well known to lie in the
  holonomy Lie algebra, this completes the proof.
\end{proof}

This evidently implies Theorem \ref{main} from the introduction. 

\subsection*{Remark} Let $M$ be a smooth manifold and let $\nabla$ be
a torsion free linear connection on $TM$ whose holonomy Lie algebra is
contained in one of the Lie algebras $\fg_0$ from Table 1. If we
further assume that there is a group $G$ such that the holonomy group
of $\nabla$ is contained in the corresponding subgroup $G_0$ (which
usually is no restriction), then $\nabla$ is induced from a first
order $G_0$--structure $E\to M$. It is a general fact that the
harmonic torsion of such a structure which was mentioned in Section
\ref{2} can be obtain as a certain component of the torsion of any
connection on the bundle. Hence we see that $E\to M$ has to have
vanishing harmonic torsion, so $\nabla$ is a Weyl connection for an
AHS--structure.

For many choices of a $|1|$--graded Lie algebra $\fg$, vanishing of the
harmonic torsion of an AHS--structure already implies local flatness,
i.e.~local isomorphism to the homogeneous model $G/P$. Indeed, the
only instances in Table 1 for which there exist non--flat torsion free
geometries are the first and fifth lines, the second and fourth lines
for $p=2$ and the third line for $p=1$. Hence in all other cases we
are, at least up to local isomorphism, in the situation of a Weyl
connection on $G/P$.

\bibliographystyle{amsplain}

\end{document}